\documentclass{amsart}
\usepackage{amsmath,amssymb,amscd,latexsym}
\usepackage{enumerate,verbatim,calc}
\usepackage[all]{xy}
\SelectTips{eu}{}

\newcommand{\quism}{\text{quasi}\-\text{iso}\-\text{mor}\-\text{phism}}

\newcommand{\gam}[2]{{\gamma^{(#1)}_{#2}}}
\newcommand{\gamp}[2]{{\gamma'{}^{(#1)}_{#2}}}
\newcommand{\Gg}[1]{{G^{(#1)}}}

\newcommand{\Gp}[1]{{G'{}^{(#1)}}}

\newcommand{\ts}{\textstyle}

\newcommand{\sss}{\scriptscriptstyle}

\newcommand{\ges}{{\sss\geqslant}}
\newcommand{\les}{{\sss\leqslant}}

\newcommand{\bse}{{\boldsymbol e}}

\newcommand{\bsi}{{\boldsymbol i}}
\newcommand{\bsj}{{\boldsymbol j}}
\newcommand{\bsp}{{\boldsymbol p}}

\newcommand{\eps}{\epsilon}

\newcommand{\card}{\operatorname{card}}

\newcommand{\col}{\colon}

\newcommand{\dd}{\partial}
\newcommand{\depth}{\operatorname{depth}}
\newcommand{\devi}[2]{\varepsilon_{#1}(#2)}
\newcommand{\ED}[3]{{}_{#1}{\operatorname{d}}^{#2,#3}}

\newcommand{\EC}[3]{{}_{#1}\!\operatorname{E}^{#2,#3}}

\newcommand{\fm}{{\mathfrak m}}
\newcommand{\fn}{{\mathfrak n}}
\newcommand{\fp}{{\mathfrak p}}
\newcommand{\fq}{{\mathfrak q}}
\newcommand{\height}{{\operatorname{height}}}

\newcommand{\hch}[4]{\operatorname{HH}_{#1}(#3\var #2;#4)}
\newcommand{\hcoh}[4]{\operatorname{HH}^{#1}(#3\var #2;#4)}

\newcommand{\hh}[1]{\operatorname{H}(#1)}

\newcommand{\HH}[2]{\operatorname{H}_{#1}(#2)}

\newcommand{\id}{\operatorname{id}}
\newcommand{\image}{\operatorname{Im}}
\newcommand{\Ker}{\operatorname{Ker}}

\newcommand{\lra}{\longrightarrow}

\newcommand{\pd}{\operatorname{pd}}

\newcommand{\rank}{\operatorname{rank}}

\newcommand{\Spec}{\operatorname{Spec}}
\newcommand{\Supp}{\operatorname{Supp}}
\newcommand{\susp}{{\mathsf\Sigma}}

\newcommand{\Ext}[4]{\operatorname{Ext}^{#1}_{#2}(#3,#4){}}
\newcommand{\Hom}[3]{\operatorname{Hom}_{#1}(#2,#3)}
\newcommand{\Rhom}[3]{\operatorname{\mathbf{R}Hom}_{#1}(#2,#3)}
\newcommand{\dtensor}[1]{\otimes^{\mathbf{L}}_{#1}}
\newcommand{\Tor}[4]{\operatorname{Tor}_{#1}^{#2}(#3,#4){}}

\newcommand{\var}{{\hskip.7pt\vert\hskip.7pt}}

\newcommand{\vf}{{\varphi}}

\newcommand{\wh}{\widehat}
\newcommand{\wt}{\widetilde}

\newcommand{\xra}{\xrightarrow}
\newcommand{\ZZ}{\operatorname{Z}}

\newcommand{\BN}{{\mathbb N}}

\newcommand{\BZ}{{\mathbb Z}}

\theoremstyle{plain}

\newtheorem{theorem}{Theorem}[section]

\newtheorem{subtheorem}{Theorem}[subsection]
\newtheorem{sublemma}[subtheorem]{Lemma}

\newtheorem*{Corollary}{Corollary}
\newtheorem*{itheorem}{Main Theorem}

\theoremstyle{definition}

\newtheorem{example}[theorem]{Example}

\newtheorem{chunk}[theorem]{}
\newtheorem{subchunk}[subtheorem]{}

\theoremstyle{remark}

\newtheorem*{Remark}{Remark}

\numberwithin{equation}{theorem}

\begin{document}

\title[Gaps in Hochschild cohomology]
{Gaps in Hochschild cohomology
imply smoothness
for commutative algebras}
\author{Luchezar L.~Avramov}
\address{Department of Mathematics,
  University of Nebraska, Lincoln, NE 68588, U.S.A.}
\email{avramov@math.unl.edu}

\author{Srikanth Iyengar}
\address{Department of Mathematics,
University of Nebraska, Lincoln, NE 68588, U.S.A.}
\email{iyengar@math.unl.edu}

\thanks{Research partly supported by NSF grants DMS 0201904 (L.L.A) and
DMS 0442242 (S.I.)}



\date{\today}

\begin{abstract}
The paper concerns Hochschild cohomology of a commutative algebra $S$,
which is essentially of finite type over a commutative noetherian ring $K$
and projective as a $K$-module.  For a finite $S$-module $M$
it is proved that vanishing of $\hcoh{n}KSM$ in sufficiently long intervals
imply the smoothness of $S_\fq$ over $K$ for all prime ideals $\fq$
in the support of $M$.  In particular, $S$ is smooth if $\hcoh{n}KSS=0$
for $(\dim S+2)$ consecutive $n\ge0$.
 \end{abstract}

\maketitle

\section*{Introduction}

Let $K$ be a commutative noetherian ring, $S$ a commutative $K$-algebra,
and $M$ an $S$-module. We let $\hch{\ast}KSM$ and $\hcoh \ast KSM$ denote,
respectively, the Hochschild homology and the Hochschild cohomology
of the $K$-algebra $S$ with coefficients in $M$.  For each $n\in\BZ$
there are canonical homomorphisms
\begin{align*}
\lambda_n^M  &\col ({\ts\wedge}_S^n  \Omega_{S\var{K}})\otimes_S M \lra \hch{n}KSM  \\
\lambda^n_M &\col \hcoh{n}KSM \lra \Hom S{{\ts \wedge}^{n}\Omega_{S\var K}}{M}
\end{align*}
of $S$-modules, where $\Omega_{S\var K}$ is the $S$-module of $K$-linear
K\"ahler differentials of $S$.  Other concepts appearing
the next result are defined following its statement.

\begin{itheorem}
\label{ihoch} Let $K$ be a commutative noetherian ring and $S$ a
commutative $K$-algebra essentially of finite type, flat
as a $K$-module.  For a prime ideal $\fq$ in $S$ and a finite
$S$-module $M$ with $M_\fq\ne0$ the following conditions are equivalent$:$
\begin{enumerate}[\quad\rm(i)]
\item[\rm(i$_{}$)]
The $K$-algebra $S_\fq$ is smooth.
\item[\rm(ii${}_*$)]
Each map $(\lambda_n^{M})_\fq$\! is bijective and the $S_\fq$-module
$\Omega_{S_\fq\var{K}}\!$ is projective.
\item[\rm(iii$_\ast$)]
There exist non-negative integers $t,u$ of different parity satisfying
\[
\hch{t}KS{M}_\fq=0=\hch{u}KS{M}_\fq
\]
 \end{enumerate}
When the $K$-module $S$ is projective they are also equivalent to:
\begin{enumerate}[\quad\rm(i)]
\item[\rm(ii${}^*$)]
Each map $(\lambda^n_{M})_\fq$ is bijective.
\item[\rm(iii$^\ast$)]
There exist non-negative integers $t,u$ of different parity satisfying
\[
\hcoh{t+i}KS{M}_\fq=0=\hcoh{u+i}KS{M}_\fq
\quad\text{for}\quad
0\le i\le \dim_{S_\fq}M_\fq
\]
 \end{enumerate}
 \end{itheorem}

Recall that one says that $S$ is \emph{essentially of finite type}
if it is a localization of a finitely generated $K$-algebra.  A flat
$K$-algebra $S$ essentially of finite type is \emph{smooth} if the
structure map $K\to S$ has geometrically regular fibers.  Equivalently,
for every homomorphism of rings $K\to\ell$, where $\ell$ is field,
the ring $S\otimes_K\ell$ has finite global dimension.   We say that
an $S$-module $M$ is \emph{finite} if it is finitely generated, and let
$\dim_SM$ denote the Krull dimension of $M$.

The theorem incorporates several known results, discussed below.
There are two new aspects to our characterizations of smoothness: the
use of \emph{cohomology} (with a couple of exceptions, earlier results
used vanishing of homology) and the introduction of \emph{coefficients}
(all earlier results dealt with the case $M=S$).  A special case of the
theorem relates to a question of Happel \cite[(1.4)]{Hap}:
\begin{enumerate}[{\ }(1)]
\item[{}]
For a (not necessarily commutative) algebra $A$ over a field $K$,
with $\rank_KA$ finite, does $\hcoh nKAA=0$ for $n\gg 0$ imply finite
global dimension?
\end{enumerate}
The next corollary provides a strong affirmative answer in the
commutative case. This is in sharp contrast to the general situation,
where the answer is negative: see the companion paper \cite{BGS} by
Buchweitz, Green, Madsen, and Solberg.

\begin{Corollary}
Let $K$ be a field let $S$ be a commutative $K$-algebra, finite
dimensional as a $K$-vector space.  If $\hcoh nKSS=0$ for two non-negative
values of $n$ of different parity, then $S$ is a product of separable
field extensions of $K$.
 \end{Corollary}

\begin{proof}
The hypothesis $\rank_KS<\infty$ implies that $\dim S$ is $0$, and that
$S$ is smooth precisely when it is a product of finite separable field
extensions of $K$.
 \end{proof}

We place our result in the context of earlier work relating vanishing of
Hoch\-schild (co)homo\-logy and smoothness.  As always, $\Spec S$  denotes
the set of prime ideal of $S$; its subset $\Supp_S M=\{\fq\in \Spec
S\var M_\fq\ne0\}$ is the \emph{support} of $M$.

\medskip

\noindent{\bf Antecedents.}\,\ 
Let  $S$ be a $K$-algebra $S$ essentially of finite type, flat as a
$K$-module.  When citing results, a roman numeral in \emph{italic font}
font indicates the variant of the correspondingly numbered condition
in the Main Theorem, where the hypothesis is assumed to hold for $M=S$
and for all $\fq\in\Spec S$.

\medskip

\noindent\emph{The HKR Theorem}.
Hochschild, Kostant, and Rosenberg \cite{HKR} (when $K$ is a perfect
field) and Andr\'e \cite{An2} (in general) proved (i) $\implies$
({ii}${}_*$) $\&$ ({ii}${}^*$).  As $S$ is essentially of finite type,
the $S$-module $\Omega_{S|K}$ is finite, so $\wedge^n_S
\Omega_{S|K}=0$ holds for all $n\gg0$, hence one always has ({ii}${}_*$)
$\implies$ ({iii}${}_*$) and ({ii}${}^*$) $\implies$ ({iii}${}^*$).

\medskip

\noindent\emph{Homological converses to the HKR Theorem}. Andr\'e
\cite{An2} proved (\emph{ii}${}_*$) $\implies$ (\emph{i}).

(\emph{iii}${}_*$) $\implies$ (\emph{i}) was proved by Avramov and
Vigu\'e-Poirrier \cite{AV} when $K$ is a field; by Campillo,
Guccione, Guccione, Redondo, Solotar, and Villamayor \cite{BACH}
when, in addition, $\operatorname{char}(K)=0$; by Rodicio
\cite{Ro3} in general.

\medskip

\noindent\emph{Cohomological converses to the HKR Theorem}. Assume
$S$ is projective over $K$.

For a Gorenstein ring $S$ Blanco and Majadas \cite{BM} proved
that $\hcoh nKSS=0$ for $(\dim S+2)$ consecutive values of $n\ge0$
implies $S$ is smooth over $K$; this is subsumed in the implication
({iii}${}^*$) $\implies$ ({i}) of the Main Theorem.  In joint work 
with Rodicio \cite{BMR} they showed that if $S$ is
locally complete intersection over $K$, then $\hcoh{2n}KSS=0$ or
$\Ker(\lambda^{2n}_S)=0$ for a single $n\ge0$ implies $S$ is smooth.

\medskip

\noindent{\bf Generalizations.}\,\ 
The Main Theorem is a special case of a much more
general result, Theorem \eqref{retracts:cicriteria}, concerning gaps
in $\Tor*RSM$ and $\Ext*RSN$ when $R$ is a noetherian ring, $S$ is an
algebra retract of $R$, and $M$ is a \emph{complex} of $S$-modules.
For $\Tor*RSS$ that result is due to Rodicio \cite{Ro3}.  However, 
to prove (iii${}^*$) $\implies$ (i) in the Main Theorem, even for $M=S$,
we do need to use complexes.

\medskip

\noindent{\bf Conventions.}\,\ 
In the rest of this article all rings are assumed to be commutative.
A \emph{local ring} is a noetherian ring that has a unique maximal ideal.
A \emph{local homomorphism} is a homomorphism of rings, whose source
and target are local and which maps maximal ideal into maximal ideal.

\section{Closed homomorphisms}
\label{Closed homomorphisms} In this section $\vf\col(R,\fm,k)\to
S$ is a surjective local homomorphism.

We recall a general construction due to Tate \cite{Ta}.  More details
about Tate resolutions and acyclic closures can be found in the original
paper, in the book of Gulliksen and Levin \cite[Chapter I]{GL}, or in
the survey \cite[Chapter 6]{barca}.

\begin{subsection}{Tate resolutions.}
\label{tate:resolutions} 
For each positive integer $n$ let $X_n$ denote a free graded $R$-module
concentrated in degree $n$; furthermore, $R\langle X_n\rangle$ denotes
the exterior algebra on $X_n$ if $n$ is odd, and the divided powers
algebra on $X_n$ if $n$ is even; in the latter case, $x^{(i)}$ denotes
the $i$th divided power of $x\in X_n$.

A \emph{Tate resolution} of $\vf$ is a DG (= differential graded) algebra
$G$ having a system of divided powers compatible with the action of
the differential and a filtration $\{\Gg p\}_{p\ge0}$ by DG subalgebras 
with divided powers, such that
\begin{enumerate}
\item[\rm(0)]
$\Gg 0=R$ and $\Gg{p-1}\subseteq\Gg p$, for $p\ge1$.
\item
$\Gg p=\Gg{p-1}\otimes_R R\langle X_p\rangle$ as graded
$R$-modules, for $p\ge1$.
\item
$\dd(x^{(i)})=\dd(x)x^{(i-1)}$ for all $i\ge1$ when $|x|$ is even and positive.
\item
$\HH{0}{\Gg p}=S$ for $p\ge1$.
\item
$\HH{i}{\Gg p}=0$ for $1\le i<p$.
\item
$G=\bigcup_{p\ges0}\Gg p$.
 \end{enumerate}
Forgetting the multiplicative structures, $G$ is a free resolution of
$R$ over $S$.  One always exists: form DG algebras satisfying
conditions (0) through (4) by induction on $p$, then use (5) to
define $G$.  Control may be exercised at each step of the
process.

As starting point, one may choose any surjective $R$-linear map
\[
\delta_1\col X_1\lra \Ker(\vf)
\]
and define the differential on $\Gg 1$ so that its restriction to $X_1$
is the composition of $\delta_1$ with the inclusion $\Ker(\vf)\subseteq
R=G^{(0)}_0$. If $\bse_1$ is a basis for $X_1$, then $\delta_1(\bse_1)$
generates the ideal $\Ker(\vf)$ and $\Gg 1$ is the Koszul complex on
$\delta(\bse_1)$.

For each $p\ge2$ one may choose any surjective $R$-linear map
\[
\delta_p\col X_p\to\HH{p-1}{\Gg {p-1}}
\]
lift it to a homomorphism $\wt\delta_p\col X_p\to
\ZZ_{p-1}(\Gg{p-1})$, and define a differential on $\Gg p$, which
on $X_p$ is the composition of $\wt\delta_p$ with
$\ZZ_{p-1}(\Gg{p-1})\subseteq G^{(p-1)}_{p-1}=G^{(p)}_{p-1}$.
 \end{subsection}

\begin{subsection}{Acyclic closures.}
\label{Acyclic closures}
An \emph{acyclic closure} of $\vf$ is a Tate resolution obtained by
choosing for each $p\ge1$ the map $\delta_p$ in \eqref{tate:resolutions}
to be a projective cover.

Let $G$ be an acyclic closure of $\vf$ and let $G'$ be a Tate resolution
of $\vf$.  There exists then a morphism $\gamma\col G\to G'$ of DG $R$-algebras
with divided powers, and for any such morphism the homomorphism of
$R$-modules $\gamma_n\col G_n\to G'_n$ is a split injection.  If $G'$
is also an acyclic closure of $\vf$, then $\gamma$ is an isomorphism,
and it induces an isomorphism $\Gg p\to \Gp p$ for each $p\ge0$.

In particular, the $p$th \emph{stage} $\Gg p$ of an acyclic closure $G$ of
$\vf$ is independent, up to isomorphism, of the choice of $G$.

The next remark is immediate from the construction of acyclic closures.

\begin{subchunk}
\label{g1=koszul} $\Gg 1$ is the Koszul complex $E$ on a minimal
generating set for $\Ker(\vf)$.
\end{subchunk}

We introduce two numerical invariants of $\vf$, for use
throughout the paper. Letting $\nu_S(N)$ denote the minimal number
of generators an $S$-module $N$, we set
\[
\devi 2\vf = \nu_S(\Ker(\vf)) \quad\text{and}\quad \devi 3\vf=\nu_S(\HH 1E)
\]
These are part of the \emph{deviations} of $\vf$; see
\cite[(2.5)]{AI:ens}. The first assertion below is clear; the second one is
a standard characterization of regular sequences.

\begin{subchunk}
$\devi 2\vf=0$ if and only if $\vf=\id^R$.
\end{subchunk}

\begin{subchunk}
\label{e3:vanishing}
$\devi 3\vf=0$ if and only if $\vf$ is generated by a regular sequence.
\end{subchunk}
 \end{subsection}
 
A complex $F$ of finite free $R$-modules is said to be \emph{minimal}
if $\dd(F)\subseteq\fm F$.  For each integer $p\ge1$, the construction
of the $p$th stage $F_{\les p}$ of a minimal resolution of $S$ adds to
$F_{\les p-1}$ a single new free module in degree $p$.

In contrast, the construction of the $p$th stage $\Gg p$ of an acyclic
closure of $\vf$ adds to $\Gg{p-1}$ shifts of \emph{every} free module
present in it: finitely many shifts appear when $p$ is odd, and infinitely
many when $p$ is even.  Thus, when the resolution of $S$ over $R$ provided
by an acyclic closure is minimal, one has a certain control of the growth
of that resolution.

This explains our interest in the class of maps described below.

\begin{subsection}{Closed homomorphisms.}
\label{closure} We say that the homomorphism $\vf$ is
\emph{closed} if some acyclic closure $G$ of $\vf$ is a minimal
resolution of $S$ over $R$.

A celebrated result of Gulliksen \cite{Gu1} and Schoeller \cite{Sc}
can be read as follows:

\begin{subchunk}
\label{residue field} The canonical surjection $R\to k$ is closed
for every $R$.
 \end{subchunk}

To state an extension, we recall that the homomorphism $\vf$ is
\emph{large} if the map
\[
\Tor{n}\vf kk \col \Tor{n}Rkk \to \Tor{n}Skk
\]
is surjective for each $n$.  The notion was introduced by Levin \cite{Le}.
The following theorem of Avramov and Rahbar-Rochandel, see \cite[(2.5)]{Le}, 
provides a significant supply of closed homomorphisms.

\begin{subchunk}
\label{large is closed} Every large homomorphism is closed.
 \end{subchunk}

The last result will be applied through the following observation:

\begin{subchunk}
\label{retracts}
If there is a homomorphism of rings $\psi\col S\to R$ with
$\vf\circ\psi=\id^S$, then
\[
\Tor{n}\vf kk \circ \Tor n\psi kk
 = \Tor n{\vf\circ\psi}kk = \Tor n{\id^S}kk = \id^{\Tor n Skk}
\]
by functoriality; thus, $\Tor n\vf kk$ is surjective, hence
$\vf$ is large, and so closed.
\end{subchunk}
  \end{subsection}

In this paper we are mostly interested in obtaining lower bounds on the
sizes of the $S$-modules $\Tor nRSM$ and $\Ext nRSM$.  For that purpose
we use properties of $\vf$ that are weaker than closure.

\begin{subsection}{Partly closed homomorphisms.}
\label{$p$-closed}
Let $\Gg p$ be as in \eqref{Acyclic closures} for some $p\ge1$ and $F$
be a minimal free resolution of the $R$-module $S$. As $\HH 0{\Gg p}=S$
and each $G^{(p)}_n$ is $R$-free, the augmentation $\Gg p\to S$
lifts to a \emph{comparison morphism}
\[
\gam p{}\col\Gg p\to F
\]
We say that $\vf$ is \emph{$p$-closed} if $\gam pn$ has an
$R$-linear left inverse for each $n\in\BZ$.

A homomorphism $\gamma$ of free $R$-modules of finite rank has a
left inverse if and only if the map $k\otimes_R\gamma$ is
injective.  This yields an alternative description:

\begin{subchunk}
\label{minimal} The homomorphism $\vf$ is $p$-closed if and only
if $\Gg p$ is minimal and the induced map $\hh{k\otimes_R\gam
p{}}\col k\otimes_R\Gg p\to \Tor{}RkS$ is injective.
 \end{subchunk}

\begin{subchunk}
\label{independence} If the homomorphism $\vf$ is $p$-closed, $G'$
an acyclic closure of $\vf$, $F'$ is a resolution of $S$ by finite
free $R$-modules, and $\gamp p{}\col\Gp p\to F'$ is a comparison
morphism, then $\gamp pn$ has a left inverse for each $n\in\BZ$.

Indeed, \eqref{Acyclic closures} yields an isomorphism $\alpha\col\Gg
p\to\Gp p$ of DG algebras over $R$, so $\Gp p$ is minimal by
\eqref{minimal}. For any comparison morphism $\beta\col F'\to F$, the
morphisms $\gam p{}$ and $\beta\gamp p{}\alpha$ are homotopic.  Thus,
$\hh{k\otimes_R\gam p{}}$ factors as
\[
k\otimes_RG^{(p)}\xra{k\otimes_R\alpha}{k\otimes_RG'^{(p)}}
\xra{\hh{k\otimes_R\gamma'^{(p)}}}\hh{k\otimes_RF'}
\xra{\hh{k\otimes_R\beta}}{k\otimes_RF}
\]
It follows that $\hh{k\otimes_R\gamma'^{(p)}}$ is injective, see
\eqref{minimal}, hence $\gamma'{}^{(p)}_n$ is split injective.
 \end{subchunk}

\begin{subchunk}
\label{base change}
Let $R'$ be a local ring, let $\rho\col R\to R'$ be a faithfully flat
local homomorphism, set $S'=R'\otimes_RS$, and let $\rho'\col R'\to S'$
denote the induced homomorphism.  The map $\vf$ is $p$-closed if and
only if so is $\vf'$.

Indeed, $R'\otimes_RG$ is an acyclic closure of $\vf'$ if
$G$ is one of $\vf$; see \cite[(1.9.8)]{GL}.
 \end{subchunk}
\end{subsection}

Next we place the properties discussed above in a familiar context,
focusing on the case $p\le2$ because these are the classes of maps
important for this paper.

\begin{subsection}{Comparisons.}
\label{ci:remarks} The homomorphism $\vf$ is said to be
\emph{complete intersection} (or \emph{c.i.}, for short) if
$\Ker(\vf)$ is generated by a regular sequence.  Clearly, one has
\begin{align*}
\text{c.i.}\implies\text{closed}
\implies\text{$2$-closed}\implies\text{$1$-closed}
\end{align*}

\begin{subchunk}
The first implication is obviously strict; see for instance \eqref{residue field}.
\end{subchunk}

\begin{subchunk}
The canonical map from $R=k[[x,y]]/(x^2,xy)$ to $S=R/(y^2)$ is
$1$-closed, but not $2$-closed: apply \eqref{independence} to the
commutative diagram
\[
\xymatrixrowsep{3pc}
\xymatrixcolsep{3pc}
\xymatrix{
\cdots\ar@{->}[r]^-{x} &R\ar@{->}[r]^-{y^2}
\ar@{->}[d]_-{\left[\begin{smallmatrix}0\\y\end{smallmatrix}\right]}
&R\ar@{->}[r]^-{x} \ar@{=}[d] &R\ar@{->}[r]^-{y^2} \ar@{=}[d]
&R\ar@{->}[r] \ar@{=}[d] &0
\\
\cdots\ar@{->}[r]
&R^2\ar@{->}[r]^-{[x\ y]}
&R\ar@{->}[r]^-{x}
&R\ar@{->}[r]^-{y^2}
&R\ar@{->}[r]
&0
}
\]
whose top row is the beginning of the second stage of an acyclic closure
of $\vf$ and whose bottom row is the beginning of a minimal resolution
of $S$ over $R$.
 \end{subchunk}

\begin{subchunk}
We do not whether every $2$-closed homomorphisms is actually closed.
\end{subchunk}
 \end{subsection}

Except for the name, $1$-closed homomorphisms have appeared in
literature.

\subsection{One-closed homomorphisms}
As defined, $p$-closure requires $\gam pn$ to be split
injective for each $n$. However, $1$-closure can be detected in a
single degree.

\begin{sublemma}
\label{sop} Let $E$ be the Koszul complex on a minimal generating
set for the ideal $\Ker(\vf)$ and let $\gamma\col E\to F$ be a comparison
morphism to a minimal free resolution of $S$ over $R$. The
following conditions are equivalent:
\begin{enumerate}[\quad\rm(1)]
\item The homomorphism $\vf$ is $1$-closed.
\item
 For $c=\devi2\vf$ the map $(k\otimes_R\gamma_c)\col(k\otimes_RE_c)\to\Tor
cRkS$ is injective.
\end{enumerate}
\end{sublemma}

\begin{proof}
Observation \eqref{g1=koszul} and property \eqref{minimal} show that (1)
implies (2). 

For the
converse, note that the isomorphism $k\otimes_R E\cong {\ts\wedge k^c}$ 
of graded $k$-algebras shows that the socle
of $k\otimes_R E$ is $k\otimes_RE_c$.  As $k\otimes_R\gamma$ is a 
homomorphism of graded $k$-algebras, when
$k\otimes_R\gamma_c$ is injective so is $k\otimes_R\gamma$; now use
\eqref{independence}.
 \end{proof}

The preceding description brings to light a connection between
$1$-closure for parameter ideals and Hochster's Canonical Element
Conjecture, see \cite{Ho}.

\begin{subchunk}
\label{cec}
Let $R$ be a local ring. The following conditions are equivalent:
\begin{enumerate}[\quad\rm(1)]
\item The Canonical Element Conjecture holds for $R$. \item For
each system of parameters $\bsp$ of $R$ the map $\vf\col R\to
R/(\bsp)$ is $1$-closed.
\end{enumerate}

Indeed, Roberts \cite{Pr} has proved that the Canonical Element
Conjecture holds for $R$ if and only if for each free resolution
$F$ of $R/(\bsp)$ over $R$ and each comparison morphism
$\kappa\col E\to F$, the induced map $(k\otimes_RE_c)\to
\HH{c}{k\otimes_RF}$ is injective; another proof of his result is
given by Huneke and Koh \cite[(1.3)]{HK}. Thus, the desired
equivalence is contained in \eqref{independence} and Lemma
\eqref{sop}.
 \end{subchunk}

The hypothesis on $R$ in the next theorem reflects the use in its
proof of a result of Bruns \cite{Br}, which in turn relies on the
Improved New Intersection Theorem.

\begin{subtheorem}
\label{ci:koszul} Let $\vf\col R\to S$ be a $1$-closed
homomorphism and assume $R$ contains a field as a subring.  If
$\pd_RS$ is finite, then $\vf$ is complete intersection.
 \end{subtheorem}

\begin{proof}
Set $I=\Ker\vf$ and $c=\devi 2\vf$.  The Koszul complex $E=\Gg1$
on a minimal generating set of $I$ yields an injection $\kappa\col
k\otimes_RE\to\Tor{}RkS$, see \eqref{minimal}.

Since $\pd_RS$ is finite, \cite[Lemma 2]{Br} yields $\kappa_i=0$
for $i>\height I$.  By construction one has $\rank_RE_1=c$, and this
implies $c\leq\height I$.  The reverse inequality always holds, due to
the Principal Ideal Theorem, hence one gets $\height I=c$.  On the
other hand, $\pd_RS<\infty$ implies that $\height I$ equals the maximal
length of an $R$-regular sequence in $I$, see \cite[(2.5)]{AF:CM}.
We concude that $I$ can be generated by an $R$-regular sequence,
as desired.
 \end{proof}

\section{Bounds on homology}
\label{Homology}

The main result of this section is a condition for a $2$-closed homomorphism to be c.i.
When $\vf$ admits a section and $M=S$ it specializes to a result of Rodicio, \cite[Theorem
1]{Ro3}.  The reason for dealing with complexes, rather than just with modules, will
become apparent in the proof of Theorem \eqref{ci:extcriteria}.

\begin{theorem} \label{ci:criteria} Let
$\vf\col R\to S$ be a $2$-closed local homomorphism and $M$ a
complex of $S$-modules with $\hh M$ degreewise finite and bounded
below.

If there exist integers $t,u\geq\inf\hh M$ of different parity, such that
\[
\Tor tRSM = 0 = \Tor uRSM
\]
then the homomorphism $\vf$ is complete intersection.
\end{theorem}

We comment on notions and notation appearing in the theorem and
its proof.

\begin{chunk}
\label{truncations} For definitions of Tor and Ext
for complexes we refer to \cite{Wi}. When their arguments are
modules (modules are always identified with complexes concentrated
in degree $0$) these are the classical derived functors. We set
\begin{align*}
\inf\hh M &= \inf\{n\var \HH nM\ne0\}\\
\sup\hh M &= \sup\{n\var \HH nM\ne0\}
\end{align*}
When $\inf \hh M$ (respectively, $\sup\hh M$) is finite we say that
$\hh M$ is \emph{bounded below} (respectively, \emph{above}). If $M$
is bounded on either side, then $\hh M\ne 0$, because $\hh M=0$ is
equivalent to $\inf \hh M=\infty$, and also to $\sup\hh M=-\infty$.

For each integer $j$ a complex $\susp^jM$ is defined by
\[
\susp^j(M)_n= M_{n-j} \qquad\text{and}\qquad \dd^{\susp^j M}_n=
(-1)^{|j|}\dd^M_{n-j}
\]

Morphisms of complexes are chain maps of degree $0$. A
\emph{quasiisomorphism} is a morphism that induces isomorphisms in
homology in all degrees; we tag quasiisomorphisms with the symbol
$\simeq$, and isomorphisms with $\cong$.
 \end{chunk}

We deduce Theorem \eqref{ci:criteria} from the following, much stronger, result.

\begin{theorem}\label{non-van:theorem}
Let $\vf\col R\to S$ be a $2$-closed local homomorphism, set $c=\devi 2\vf$ and
$d=\devi 3\vf$. If $M$ is a complex of $S$-modules with $\hh M$ degreewise finite and
bounded below, then  for $i=\inf \hh M$ and $m=\nu_S(\HH iM)$ one has inequalities
 \begin{alignat*}{2}
\nu_S\big(\Tor{n+i}RSM\big) &\geq m\cdot\binom cn& &\text{for}\quad 0\leq
n\leq c
\\
\nu_S\big(\Tor{2n+i+c}RSM\big) &\geq m\cdot \binom{n+d-1}{d-1}
&\quad&\text{for}\quad 1\leq n
 \end{alignat*}
  \end{theorem}

The proof uses a general lemma in homological algebra, presented
below.

\setcounter{subsection}{3}

\begin{subsection}{}
\label{functor} Let $T$ be a covariant additive functor from the
category of complexes of $S$-modules to the category of graded
$S$-modules; for each complex $M$ of $S$-modules we write $T_n(M)$
for the component in degree $n$ of the graded $S$-module $T(M)$.
Assume, furthermore, that $T$ has the following properties:
\begin{enumerate}[\quad\rm(a)]
\item 
$T$ preserves quasiisomorphisms.
\item 
$T$ commutes with shifts: $T_n(\susp^jM)=T_{n-j}(M)$ for each $n\in\BZ$.
\item 
$T(\mu^M_s)=\mu_s^{T(M)}$ for each $s\in S$, where $\mu_s$ denotes
multiplication by $s$.
 \end{enumerate}

For the maximal ideal $\fn$ of $S$ property (c) implies:

\begin{subchunk}
One has $\fn\cdot T(k)=0$, so $T(k)$ is naturally a graded
$k$-vector space.
\end{subchunk}

\begin{sublemma}
\label{functor:lemma}
Let $\vf\col R\to S$ be a surjective homomorphism and $\eps\col S\to k$
the canonical surjection.  If $M$ is a complex of $S$-modules
as in Theorem \eqref{non-van:theorem}, then
\[
\nu_S(T_{n+i}(M))\geq m\cdot \rank_k \image(T_n(\eps))
\]
\end{sublemma}

\begin{proof}
 First we simplify $M$. The inclusion into $M$ of the subcomplex
\[
M':=\quad \cdots \lra M_{i+2} \xrightarrow{\dd_{i+2}} M_{i+1}
\lra \Ker(\dd_i) \lra 0
\]
is a quasiisomorphism. By (a) one has $T_n(M) {\cong} T_n(M')$,
so we may assume $M=M'$.  Set $H=\HH iM$, choose a surjection $H\to k^m$
and let $\pi$ denote the composition $M \to \susp^i H \to \susp^i k^m$ of
morphisms of complexes of $S$-modules.  Lifting $\susp^i \eps^m$ over $\pi$
to a morphism $\rho\col \susp^iS^m\to M$, we get a commutative diagram
\[ \xymatrixrowsep{1.5pc}
\xymatrixcolsep{3pc} 
\xymatrix{ 
& T_{n+i}(\susp^i S^m) \ar@{->}[dl]_-{T_{n+i}(\rho)}\ar@{->}[r]^{\cong}
  &T_{n+i}(\susp^iS)^{m}\ar@{->}[dd]^-{T_{n+i}(\susp^i\eps)^m} \\
T_{n+i}(M) \ar@{->}[dr]_-{T_{n+i}(\pi)}\\
  &T_{n+i}(\susp^ik^m) \ar@{<-}[uu]_-{T_{n+i}(\susp^i\eps^m)}\ar@{->}^-{\cong}[r]
  &T_{n+i}(\susp^ik)^{m}
}
\]
of homomorphism of $S$-modules. We can now write the relations below
\begin{align*}
  \nu_S(T_{n+i}(M))
  &\ge \rank_k\image(T_{n+i}(\pi))\\
  &\geq \rank_k \image(T_{n+i}(\susp^i\eps^m)) \\
  &= m\cdot \rank _k\image(T_{n+i}(\susp^i\eps))\\
  &= m\cdot \rank _k\image(T_n(\eps))
\end{align*}
by using consecutively the following facts: the maximal ideal of $S$ annihilates
$T_{n+i}(\susp^ik)$; the diagram commutes; $T_{n+i}(\susp^ik)$ is isomorphic to
$T_{n}(k)$.
\end{proof}
\end{subsection}

We need an explicit description of a subcomplex of an
acyclic closure $G$ of $\vf$.

\setcounter{theorem}{4}

\begin{chunk}
\label{tate complex}
In the notation of \eqref{tate:resolutions}, each $R$-module $G^{(2)}_n$ has a basis
\[
\bigg\{x_\bsi y_\bsj\ \bigg|
\ \bsi\subseteq[1,c]\,,\bsj=(j_1,\dots,j_d)\in
\BN^d\,,\ \card \bsi+2\sum_{h=1}^d j_h=n\bigg\}
\]
where $[1,c]=\{1,\dots,c\}$.  Let $a_1,\dots,a_c$ be a
minimal set of generators of $\Ker\vf$ and $E$ the Koszul complex on
it.  The differential of $\Gg 2$ then has the form
\[
\dd\big(x_{\bsi}y_{\bsj}\big)= \sum_{i\in \bsi} \pm a_{i}\,
x_{\bsi\setminus\{i\}}y_\bsj +\sum_{i\in[1,c]}\,\sum_{j=1}^d\pm
b_{ij}\,x_{\bsi\cup\{i\}}y_{\bsj-\bse_j}
\]
where $\bse_j\in\BN^d$ is the $j$th unit vector, and
\[
z_j=\sum_{i=1}^c b_{ij}\,x_{i}\in G^{(2)}_1\qquad\text{for}\qquad
j=1,\dots,c
\]
are cycles whose homology classes minimally generate $\HH1{E}$.

All the coefficients $a_i$ and $b_{ij}$ are in $\fm$: this is clear
for the $a_i$; as they form a minimal set of generators the relation
$0=\dd(z_j)=\sum_{i=1}^c b_{ij}\,a_{i}$ implies $b_{ij}\in\fm$.
 \end{chunk}

\begin{proof}[Proof of Theorem \emph{\eqref{non-van:theorem}}]
Let $\eps\col S\to k$ be the canonical surjection. Lemma
\eqref{functor:lemma} applied to the functor $T$ defined by
$T_n(M)=\Tor nRSM$ yields
\[
\nu_S(\Tor{n+i}RSM)\ge m\cdot \rank_k\image(\Tor nRS\eps)
\]
Next we estimate the rank on the right hand side. Let $\Gg 2$ be
the second stage in an acyclic closure of $\vf$, $F$ a free
resolution of $S$ over $R$, and let $\gamma\col \Gg 2\to F$ be a
comparison morphism. The following diagram commmutes:
\[
\xymatrixrowsep{3pc}
\xymatrixcolsep{2pc} 
\xymatrix{ 
\HH n{\Gg 2\otimes_RS} \ar@{->}[d]_{\HH n{\Gg 2\otimes_R{\eps}}}
\ar@{->}[rr]^{\HH n{\gamma\otimes_RS}}
   && \HH n{F\otimes_RS} \ar@{=}[r] \ar@{->}[d]^{\HH n{F\otimes_R\eps}}
     & \Tor{n}RSS \ar@{->}[d]^{\Tor{n}RS\eps}  \\
  \HH n{\Gg 2\otimes_Rk} \ar@{->}[rr]^{\HH n{\gamma\otimes_Rk}} && \HH n{F\otimes_Rk}
  \ar@{=}[r] & \Tor{n}RSk }
\]
As $\HH n{\gamma\otimes_Rk}$ is injective by \eqref{minimal}, for
each $n$ we get
\[
\rank_k(\image\Tor{n}RS{\eps})\geq \rank_k \image(\HH n{\Gg 2\otimes_R\eps})
\]

{}From the description of $\Gg 2$ in \eqref{tate complex} one sees
that the graded submodule
\[
Z=\bigoplus_{\bsi\subseteq[1,c]} S(x_{\bsi}y_0\otimes1)
\oplus\bigoplus_{\bsj\in\BN^d\setminus0}
S(x_{[1,c]}y_\bsj\otimes1)\subseteq \Gg 2\otimes_RS
\]
consists of cycles and the differential of $\Gg 2\otimes_Rk$ is trivial; thus
the composition
\[
Z\otimes_Sk\lra k\otimes_S\hh{\Gg 2\otimes_RS}\lra\hh{\Gg 2\otimes_Rk}=\Gg 2\otimes_Rk
\]
is injective.  Counting ranks over $k$ one obtains inequalities
\begin{alignat*}{2}
  \rank_k\image(\HH n{\Gg 2\otimes_R\eps})&\geq\binom cn&  &\text{for}\quad 0\leq n\leq c\\
  \rank_k\image(\HH{2n+c}{\Gg 2\otimes_R\eps}) &\geq\binom{n+d-1}{n} &\quad &\text{for}\quad
  1\leq n
\end{alignat*}

To get the desired result, concatenate the (in)equalities established above.
 \end{proof}

\begin{proof}[Proof of Theorem \emph{\eqref{ci:criteria}}]
\label{ci:criteria:proof}
By hypothesis, one has $\Tor tRSM=0=\Tor uRSM$ for integers $t,u$
satisfying $t,u\ge\inf(\hh M)=i>-\infty$ and $t\not\equiv u\pmod 2$. The
first inequality established in Theorem \eqref{non-van:theorem} implies
$t,u>i+c$ for $c=\devi 2\vf$.  For $d=\devi 3\vf$ its second inequality
in the theorem then yields $\binom {n+d-1}{d-1}=0$ for some $n\ge1$,
forcing $d=0$. Thus, $\vf$ is complete intersection by \eqref{e3:vanishing}.
 \end{proof}

\section{Vanishing of cohomology}
\label{Cohomology}

In this section we provide cohomological criteria for a $2$-closed
homomorphism to be c.i.\ This uses a notion of \emph{depth} of a
complex $M$, defined by
\[
\depth_SM = \inf\{n\in\BZ \mid \Ext nSkM\ne 0\}
\]
This is the classical concept when $M$ is a finite $S$-module.  

\begin{theorem}
\label{ci:extcriteria} Let $\vf\col R\to S$ be a $2$-closed
homomorphism and $M$ a complex of $S$-modules with $\hh M$
degreewise finite and bounded above.

If there exist integers $t,u\geq \depth_S M - \dim S$, of different
parity, such that
\[
\Ext {t+n}RS M = 0 = \Ext {u+n}RSM
\quad\text{for}\quad 0\leq n\leq  \max\{\dim_S \HH nM\mid n\in\BZ\}
\]
then the homomorphism $\vf$ is complete intersection.
\end{theorem}

\begin{Remark}
As one always has $\dim S-\depth_S M\ge\sup\hh M$, see \cite[(2.11.3)]{FI},
the bound on $t,u$ in the theorem may be replaced by $t,u\ge-\sup\hh M$.
 \end{Remark}

Theorem \eqref{ci:extcriteria} is a cohomological counterpart to Theorem
\eqref{ci:criteria}, which provides a main ingredient in its
proof.  Another component is the use of properties of dualizing
complexes, reviewed below; we refer to Hartshorne \cite{Ha} for details.

\setcounter{subsection}{1}

\subsection{Dualizing complexes}
\label{dualizing}
A \emph{dualizing complex} for $(S,\fn,k)$ is a complex
\[
D=\  0\to D_0 \to D_{-1}\to \cdots \to D_{-\dim S}\to 0
\]
of injective modules with $\hh D$ degreewise finite and
$\Hom SkD\simeq\susp^{-\dim S}k$.

Up to a \quism\, of complexes, $S$ has at most one dualizing complex.
Such a complex exists when the local ring $S$ is complete.

For each complex of $S$-modules $M$ we set $M^\dagger=\Hom SMD$.

\begin{subchunk}
If $\hh M$ is degreewise finite, then so is $\hh{M^\dagger}$.
\end{subchunk}

\begin{subchunk}
If $\hh M$ is bounded on one side, then $\hh{M^\dagger}$ is bounded on the other.
\end{subchunk}

\begin{sublemma}
\label{dagger inf}
If $\hh M$ is degreewise finite and bounded above, then
\[
\inf\hh{M^\dagger}=\depth_SM-\dim S
\]
\end{sublemma}

\noindent\emph{Proof}.
The complex $\hh{M^\dagger}$ is degreewise finite and bounded below,
see \eqref{dualizing}.  This implies, the first equality
below; the second one holds by definition:
\[
\inf\hh{M^\dagger}=\inf\hh{k\dtensor S{M^\dagger}}=\inf\hh{k\dtensor S\Hom SMD}
\]

To compute the right hand side we use a sequence of quasiisomorphisms:
\begin{align*}
k\dtensor S \Hom SMD
&\simeq\Hom S{\Rhom SkM}D\\
&\simeq\Hom S{\Ext{} SkM}D\\
&\cong \Hom k{\Ext{}SkM}{\Hom SkD}\\
&\simeq\Hom k{\Ext{} SkM}{\susp^{-\dim S}k}\\
&\simeq\susp^{-\dim S}\Hom k{\Ext{} SkM}k
\end{align*}
The first one holds because $k$ has a resolution by finite free
$S$-modules and $D$ is a bounded complex of injectives.  For the second,
note that $\Rhom SkM$ can be represented by a complex of $S$-modules
annihilated by $\fn$, so it is quasiisomorphic to its own homology, namely,
$\Ext{} SkM$.  The third one holds because $\Ext{}
SkM$ is a direct sum of copies of shifts of $k$.  The fourth quasiisomorphism
is induced by $\Hom SkD\simeq\susp^{-\dim S}k$; see \eqref{dualizing}.
The last one is standard.

We now finish the computation of $\inf\hh{M^\dagger}$ as follows:
\begin{xxalignat}{3}
&{\phantom{\square}}
&\inf\hh{M^\dagger}
&=\inf\left(\susp^{-\dim S}\Hom k{\Ext{} SkM}k\right)
&&{\phantom{\square}}\\
&{\phantom{\square}}
&&=\inf\Hom k{\Ext{}SkM}k -\dim S
&&{\phantom{\square}}\\
&{\phantom{\square}}
&&=\depth_SM-\dim S
&&{\square}
 \end{xxalignat}

\begin{subchunk}
\label{supp:dualizing}
For every finite $S$-module $N$ one has
\[
\Ext nSND=0\,, \quad\text{unless}\quad \dim S - \dim_SN\le n\le \dim S-\depth_SN
\]
\end{subchunk}

\setcounter{theorem}{2}

\begin{chunk}
\label{supports}
The \emph{support} of a complex $M$ is defined to be the set
\[
\Supp_S M = \{\fq\in \Spec S \mid \hh{M_\fq}=0\}
\]
Let $\dim \Supp_SM$ denote the dimension of space $\Supp_SM$ in the Zariski
topology on $\Spec S$. It is not hard to see that if $\hh M$ is degreewise finite, then
\[
\dim\Supp_SM = \max\{\dim_S \HH nM\mid n\in\BZ\}
\]
\end{chunk}

\begin{proof}[Proof of Theorem \emph{\eqref{ci:extcriteria}}]
By \eqref{base change}, the map $\wh\vf\col\wh R\to\wh S$ of
maximal-ideal-adic completions induced by $\vf$ is $2$-closed. For
$\wh M=M\otimes_S\wh S$ and each $n\in\BZ$ one has
\[
\HH n{\wh M}\cong\HH nM\otimes_R\wh R \qquad\text{and}\qquad \Ext
{n}{\wh R}{\wh S}{\wh  M}\cong\Ext {n}RSM \otimes_R \wh R
\]
where the first one is due to the flatness of $\wh R$ over $R$, while the second uses, in
addition, that $S$ is finite over $R$ and that $\hh M$ is bounded above.   In particular,
one has $\dim_{\wh S}\HH n{\wh M}=\dim_S\HH n{M}$ for each $n\in\BZ$.  Thus, the 
hypotheses of the theorem do not change when $R$, $S$, $M$ are replaced by
$\wh R$, $\wh S$, $\wh M$, respectively.  Furthermore, if $\wh\vf$ is c.i., then so is
$\vf$.  Thus, we may assume that the ring $S$ is complete, and hence that it has a dualizing
complex $D$.  Set $m= \max\{\dim_S \HH nM\mid n\in\BZ\}$.

As $D$ is a bounded complex of injectives, there is a natural quasiisomorphism
\begin{equation}
\label{natural}
\tag{$*$}
\Hom S{\Rhom RSM}D \simeq S\dtensor R \Hom SMD
\end{equation}
The composition of the factors on the left gives rise to a spectral sequence with
\begin{gather*}
  \EC 2pq = \Ext {-p}S{\Ext {q}RSM}D\qquad\text{and}\qquad \ED rpq\col \EC rpq \lra \EC
  r{p-r}{q+r-1}
\end{gather*}
As the $R$-module $S$ is finite, one has $\Ext qRSM_\fq\cong \Ext qRS{M_\fq}$ for each
$\fq\in\Spec S$, so $\Supp_S\Ext qRSM\subseteq \Supp_S M$, for each $q\in\BZ$.  Thus, one gets
\[
\dim_{S_\fq}(\Ext qRSM_\fq)\le\dim\Supp_S M=m
\]
where the equality comes from \eqref{supports}.  Now \eqref{supp:dualizing} yields
\[
\EC 2pq=0 \quad\text{for}\quad p\notin [-\dim S, -\dim S + m]
\]
so the sequence converges.  Formula  \eqref{natural} shows that its abutment is equal to
\[
\hh{S\dtensor R \Hom SMD}=\Tor {p+q}RS{M^\dagger}
\]

On the other hand, our hypothesis entails $\EC 2pq = 0$ for
\[
t\leq q\leq t+m \quad\text{and}\quad u \leq q\leq u +m
\]
As a consequence, one obtains equalities
\[
\EC 2pq=0 \quad\text{whenever}\quad p+q= t \quad\text{or}\quad p+q= u
\]
They imply $\EC{\infty}pq=0$ if $p+q=t$ or $p+q=u$, so convergence yields
\[
\Tor {t}RS{M^\dagger}=0=\Tor {u}RS{M^\dagger}
\]

In view of Lemma \eqref{dagger inf} and our hypothesis, the
complex $M^\dagger$ satisfies
\[
\inf\hh{M^\dagger}=\depth_SM^\dagger-\dim S\le\min\{t,u\}
\]
Now Theorem \eqref{ci:criteria}, applied to $M^\dagger$, shows that $\vf$
is complete intersection.
\end{proof}

\section{(Co)homology of algebra retracts}
\label{Algebra retracts} Let $\vf\col R\to S$ be a homomorphism of
noetherian rings.

A \emph{section} of $\vf$ is a homomorphism of rings $\psi\col
S\to R$ such that $\psi\circ\vf=\id^S$; when such a homomorphism
exits $S$ is said to be an \emph{algebra retract} of $R$. Another
way to describe this situation is to say that $R$ is a
\emph{supplemented} algebra over $S$. The study of homological and
cohomological properties of supplemented algebras is a central
topic in the classical literature on homological algebra.

Each $\fq \in \Spec S$ defines a local homomorphism $\vf_\fq\col
R_{\vf^{-1}(\fq)}\to S_\fq$. If $\psi$ is a section of $\vf$, the
local homomorphism $\psi_{\fp}$, where $\fp= \vf^{-1}(\fq)$, is a
section of $\vf_\fq$. In particular, the homomorphism $\vf_\fq$ is
$2$-closed; see \eqref{retracts}.

Next we establish global versions of results from the preceding sections. To
this end we recall the construction of certain canonical homomorphisms.

\begin{chunk}
\label{kunneth} With $I=\Ker(\vf)$, one has a canonical of
$S$-modules isomorphism
\[ I/I^2\cong\Tor 1{R}SS \]
The graded $S$-module $\Tor *RSS$ has a natural structure of a strictly
commutative graded $S$-algebra, see \cite[(2.7.8)]{Wi}, so there is a
homomorphism of graded $S$-algebras: $\lambda^S\col\wedge_S(I/I^2)\to\Tor
{}RSS$.  Define $\lambda^M$ to be the composition
\[
\wedge_S(I/I^2)\otimes_S\hh M\xra{\lambda^S\otimes_S{\hh M}}\Tor{}RSS
\otimes_S\hh M \lra \Tor{}RSM
\]
where the second arrow is a K\"unneth map. Let $\lambda_M$ denote the composition
\[
\Ext {}RSM\lra\Hom S{\Tor{}RSS}{\hh M}\lra\Hom
S{{\ts\wedge}_S(I/I^2)}{\hh M}
\]
where the first arrow is a K\"unneth map, and the second is $\Hom S{\lambda^S}{\hh M}$
\end{chunk}

\begin{theorem}
\label{retracts:cicriteria}
Let  $\vf\col R\to S$ be a homomorphism of rings that admits a section, and let $M$
be a complex of $S$-modules with $\hh M$ finite. Set $I=\Ker(\vf)$.

For each prime ideal $\fq\in\Supp_SM$ the following conditions are
equivalent.
\begin{enumerate}[\quad\rm(i)]
\item
 The homomorphism $\vf_{\fq}$ is complete intersection.
\item[\rm(ii${}_*$)]
 The map $(\lambda^M)_\fq$ is bijective
 and the $S_\fq$-module $(I/I^2)_\fq$ is projective.
\item[\rm(iii${}_*$)] For integers $t,u\ge\inf\hh M_\fq$ of
different parity one has
 \[
 \Tor {t}RSM_\fq = 0 = \Tor{u}RSM_\fq
\]
\item[\rm(ii${}^*$)]The map $(\lambda_M)_\fq$ is bijective.
\item[\rm(iii${}^*$)] For integers $t,u\ge\depth_{S_\fq}M_\fq-\dim
S_\fq$ of different parity one has
 \[
 \Ext {t+i}RSM_\fq = 0 = \Ext{u+i}RSM_\fq \quad\text{for}\quad i=0,\dots,\dim\Supp_{S_\fq}M_\fq
\]
 \end{enumerate}
\end{theorem}

\begin{proof}
  Set $\fp=\vf^{-1}(\fq)$.  The $R$-module $S$ is finite, so for each $n\in\BZ$
one has
\[
\Tor {n}RSM_\fq \cong \Tor {n}{R_\fp}{S_\fq}{M_\fq}
\quad\text{and}\quad \Ext {n}RSM_\fq \cong \Ext {n}{R_\fp}{S_\fq}{M_\fq}
\]
Therefore, each of the conditions listed above
is local; moreover, any section of $\vf$ localizes to a section of
$\vf_\fq$.  Thus, changing notation we may assume that $\vf$ is a
local homomorphism and that $\fq$ is the maximal ideal of $S$.

(i) $\implies$ (ii$_{\ast}$) and (ii$^{\ast}$).  Let $E$ be the
Koszul complex on a minimal generating set of $I$.  It satisfies
$\dd(E)\subseteq IE$ and is a free resolution of $S$, as $\vf$ is
c.i.  We get
\[
\Tor{}RSS\cong \hh{E\otimes_RS}=E\otimes_RS
\]
so the $S$-module $\Tor 1RSS$ is free and $\lambda^S$ is
bijective.  Thus, $\Tor nRSS$ is finite free and vanishes for
$n<0$ or $n>\nu_S(I)$, so the K\"unneth homomorphisms
\begin{gather*}
\Tor{}RSS\otimes_S\hh M\lra \Tor{}RSM \\
\Ext {}RSM\lra\Hom S{\Tor{}RSS}{\hh M}
\end{gather*}
are bijective. The definitions of $\lambda^M$ and $\lambda_M$ show
that they are bijective as well.

(ii$_\ast$) $\implies$ (iii$_\ast$), and (ii$^\ast$) $\implies$
(iii$^\ast$). These implications are clear because the $S$-module
$I/I^2$ is finite and $\hh M$ is bounded.

(iii$_\ast$) or (iii$^\ast$) $\implies$ (i).  As $\vf$ is closed
by \eqref{retracts}, Theorem \eqref{ci:criteria}, respectively,
Theorem \eqref{ci:extcriteria}, shows that condition (iii$_\ast$),
respectively, (iii$^\ast$), implies $\vf$ is c.i.
 \end{proof}

\section{Hochschild (co)homology}

Finally, we return to the subject in the title of this article.
First we recall a classical interpretation of the functors in
question; see e.g.\ \cite[(9.1.5)]{Wi}.

\begin{chunk}
\label{hoch=tor}
Let $\eta\col K\to S$ be a homomorphism of rings and let $\vf^S\col
S\otimes_KS \to S$ be the homomorphism of rings given by $\vf^S(s\otimes_K
s')=ss'$.  

If the $K$-module $S$ is flat, then for each $n\in\BZ$ one has
\[
\hch nKSM=\Tor n{S\otimes_KS}SM
\]
Note that $I/I^2$, where $I=\Ker(\vf^S)$, is a
standard realization of the module of differentials
$\Omega_{S|K}$. The map $\lambda^M$ from \eqref{kunneth} yields
$S$-linear maps
\[
\lambda^M_n\col \wedge_S^n\Omega_{S|K}\otimes_SM\lra \hch nKSM
\]

If $S$ is projective as a $K$-module, then also
\[
\hcoh nKSM= \Ext n{S\otimes_KS}SM
\]
so in this context the homomorphism $\lambda_M$ from \eqref{kunneth} reads
\[
\lambda_M^n\col \hcoh nKSM\lra \Hom
S{{\ts\wedge}_S^n\Omega_{S|K}}M
\]
 \end{chunk}
 
 The maps above are the homomorphisms that appear in the introduction. For the proof of
 the theorem stated there we need a characterization of smoothness proved by Andr\'e
 \cite[Proposition C]{An2}, using Andr\'e-Quillen homology. A short version of his
 argument may be found in \cite[(1.1)]{AI:tata}.

\begin{chunk}
\label{smooth:diagonal}
A flat algebra $S$ essentially of finite type over a noetherian ring
$K$ is smooth if and only if the homomorphism $(\vf^S)_\fp$ is c.i.,
for each $\fp\in\Spec S$.  \end{chunk}

\begin{proof}[Proof of the Main Theorem]
Let $\vf^S\col S\otimes_KS \to S$ be the product map. We
claim that, for a given $\fq\in\Spec S$, condition (i): the
$K$-algebra $S_\fq$ is smooth, is equivalent to: (i$'$) the
homomorphism $(\vf^S)_\fq\col(S\otimes_KS)_{(\vf^S)^{-1}(\fq)}\to S_\fq$
is c.i.

Indeed, $(\vf^S)_\fq$ is surjective, so it is c.i.\ if and only if
$(\vf^S)_\fp$ is c.i.\ for each $\fp\subseteq\fq$. However, the local
homomorphisms $(\vf^S)_{\fp}$ and $(\vf^{S_\fq})_{\fp}$ coincide,
and the latter is c.i.\ for each $\fp$ precisely when the $K$-algebra
$S_\fq$ is smooth, by \eqref{smooth:diagonal}.

Given this translation and the identifications in \eqref{hoch=tor}, the
desired result is contained in Theorem \eqref{retracts:cicriteria}, for
$s\mapsto 1\otimes s$ gives a section $S\to S\otimes_KS$ of $\vf^S$.
 \end{proof}

The example below shows that condition (ii$_\ast$) in the Main Theorem
cannot be weakened in general. We do not know whether the conclusion
of the theorem still holds if the vanishing intervals in condition
(iii$^\ast$) are shortened.

\begin{example}
\label{gaps:examples}
Let $S=\BZ[\sqrt 2]$. The Hochschild homology of $S$ over $\BZ$ is
\[
\hch n{\BZ}SS =\begin{cases}S     & \text{ for  $n=0$}\cr
                              S/(2\sqrt 2) & \text{ for  odd $n\ge1$}\cr
                              0     & \text{ for  even $n\ge2$}\end{cases}
\]
while the Hochschild cohomology of $S$ over $\BZ$ is given by
\[
\hcoh n{\BZ}SS =\begin{cases} S     & \text{ for  $n=0$}\cr
                              0     & \text{ for  odd $n\ge1$}\cr
                              S/(2\sqrt 2)& \text{ for  even $n\ge2$}\end{cases}
\]

Indeed, $\Ker(S\otimes_\BZ S\to S)$ is generated by $\sqrt 2\otimes
1 - 1\otimes \sqrt 2$.  A free resolution of $S$ as a module over
$S\otimes_\BZ S$ is given by the complex $F$ below:
\[
\cdots \lra S\otimes_\BZ S
\xra{\sqrt 2\otimes 1 - 1\otimes \sqrt 2} S\otimes_\BZ S
\xra{\sqrt 2\otimes 1 + 1\otimes \sqrt 2} S\otimes_\BZ S
\xra{\sqrt 2\otimes 1 - 1\otimes \sqrt 2} S\otimes_\BZ S
\lra 0
\]
As $S$ is finite free as a $\BZ$-module, $\hch *{\BZ}SS$ is the homology of the
complex
\[
F\otimes_{S\otimes_\BZ S}S=\quad
\cdots \lra S\xra{2\sqrt 2} S \xra 0 S \xra{2\sqrt 2} S \xra {0} S\lra 0
\]
and $\hcoh *{\BZ}SS$ is the homology of the complex
\[
\Hom{S\otimes_\BZ S}FS=\quad
0\lra S \xra{0} S\xra{2\sqrt 2} S \xra{0} S \xra {2\sqrt 2} S\lra\cdots
\]
see \eqref{hoch=tor}.  The desired expressions follow.
\end{example}

\end{document}